\documentclass[leqno,draft]{article}
\usepackage{amssymb,amsmath,theorem,a4wide}

\allowdisplaybreaks[2]

\newcommand\LOR{\bigvee}
\newcommand\ET{\bigwedge}
\newcommand\model{\vDash}

\newcommand\p[1]{\langle#1\rangle}
\newcommand\nrm[1]{\lVert#1\rVert}
\newcommand\bez{\smallsetminus}
\newcommand\sset{\subseteq}
\newcommand\Sset{\supseteq}
\newcommand\nul{\varnothing}

\newcommand\num{\overline}
\newcommand\bnum{\underline}
\newcommand\eqs{\bumpeq}
\newcommand\neqs{\not\eqs}
\newcommand\red{\widetilde}
\newcommand\lead{\longrightarrow}
\newcommand\leadt{\stackrel*\lead}

\DeclareMathOperator\dom{dom}
\DeclareMathOperator\sol{Sol}
\DeclareMathOperator\des{d}

\newcommand\np{\mathrm{NP}}
\newcommand\dio{\mathrm D_}

\newcommand\M{\mathit}
\newcommand\io{\M{IOpen}}
\newcommand\idz{I\Delta_0}
\newcommand\PA{\M{PA}}
\newcommand\LQ{L_Q}

\newcommand\N{\mathbb N}
\newcommand\Z{\mathbb Z}

\newenvironment{Pf}
  {\par\noindent\textit{Proof:}\bme\ignorespaces}
  {\noproof\pagebreak[2]\vskip\medskipamount\ignorespacesafterend}
\newcommand\noproof{\leavevmode\unskip\bme\vadjust{}\nobreak\hfill$\qed$\par}
\newcommand\qed{\Box}
\def\bme{\hskip.75em\relax}

\theoremstyle{plain}
\newtheorem{Thm}{Theorem}[section]
\newtheorem{Prop}[Thm]{Proposition}
\newtheorem{Cor}[Thm]{Corollary}
\newtheorem{Lem}[Thm]{Lemma}
\newtheorem{Prob}[Thm]{Problem}
\newtheorem{Cl}{Claim}[Thm]
\renewcommand\theCl{\arabic{Cl}}

\newenvironment{Pf*}{\let\qed\qedCl\Pf}\endPf

\theorembodyfont\upshape
\newtheorem{Def}[Thm]{Definition}
\newtheorem{Exm}[Thm]{Example}

\author{Emil Je\v r\'abek\\[\medskipamount]
Institute of Mathematics of the Czech Academy of Sciences\\
\small \v Zitn\'a 25,
115\:67 Praha 1,
Czech Republic,
email: \texttt{jerabek@math.cas.cz}
}

\title{Division by zero}

\begin{document}
\maketitle

\begin{flushright}
\emph{To Albert Visser}
\end{flushright}

\begin{abstract}
For any sufficiently strong theory of arithmetic, the set of Diophantine
equations provably unsolvable in the theory is algorithmically undecidable, as a consequence of the MRDP theorem. In contrast, we show decidability of
Diophantine equations provably unsolvable in Robinson's arithmetic~$Q$. The argument hinges on an analysis of a particular
class of equations, hitherto unexplored in Diophantine literature. We also axiomatize the universal fragment of~$Q$ in
the process.


\end{abstract}

\section{Introduction}

The standard G\"odel--Church--Turing--Rosser undecidability theorem tells us that if $T$ is any consistent theory
extending Robinson's arithmetic~$Q$, the set of $\Pi_1$~consequences of~$T$ is undecidable. Furthermore, the
Matiyasevich--Robinson--Davis--Putnam theorem shows that every $\Pi_1$~formula is equivalent to unsolvability of a
certain Diophantine equation. Since the MRDP theorem can be formalized in $\idz+\M{EXP}$ due to Gaifman and
Dimitracopoulos~\cite{mrdpea}, we see that if $T$ extends $\idz+\M{EXP}$, it is undecidable whether a given Diophantine
equation is provably unsolvable in~$T$, or dually, whether it has a solution in a model of~$T$.

Surprisingly, Kaye~\cite{kaye:dio,kaye:hilb} proved that the same holds already for extensions of the weak theory~$IU_1^-$
(induction for parameter-free bounded universal formulas),
despite that it likely does \emph{not} formalize the MRDP theorem as such. One can check that Kaye's methods also
apply to extensions of Cook's theory~$\M{PV}$ of polynomial-time functions (see e.g.\ Kraj\'\i\v cek~\cite{book} for
a definition).

Going further down, decidability of solvability of Diophantine equations in models of the theory~$\io$ of quantifier-free
induction has remained an intriguing open problem ever since it was posed by Shepherdson~\cite{sheph}, see e.g.\
\cite{wil:iop,vdd:dor,otero} for partial results.

The purpose of this note is to show that solvability of Diophantine equations in models of~$Q$ is decidable,
specifically $\np$-complete. Since $Q$
does not include ring identities that allow the usual manipulations of polynomials, it may be ambiguous what exactly is
meant by Diophantine equations, so let us first state the problem precisely.
\begin{Def}\label{def:dioph}
A \emph{Diophantine equation} is a formula of the form
\[t(\vec x)=u(\vec x),\]
where $t$ and~$u$ are terms in the basic language of arithmetic $\LQ=\p{0,S,+,\cdot}$. If $T$ is a theory whose language
contains~$\LQ$, the \emph{Diophantine satisfiability problem for~$T$,} denoted $\dio T$, consists of all Diophantine
equations $t=u$ satisfiable in a model of~$T$ (shortly: $T$-satisfiable). That is,
\[\dio T=\{\p{t,u}:T+\exists\vec x\,t(\vec x)=u(\vec x)\text{ is consistent}\}.\]
\end{Def}

The decidability of~$\dio Q$ is on the whole not so surprising, as $Q$ has models with ``black holes'' (to use Albert
Visser's term) that can serve to equate nearly any pair of terms. It turns out however that while this argument yields a simple proof of decidability of
Diophantine satisfiability for certain mild extensions of~$Q$, it does not suffice for $Q$~itself: it only provides a
reduction to systems of equations of a special form (see Eq.~\eqref{eq:1} below) that we have to investigate in detail.

We need to embed certain models in models of~$Q$ as a part of our main construction, and to facilitate
this goal, we will explicitly axiomatize the universal consequences of~$Q$, which could be of independent
interest.

See \cite{doy-con,doy-qiu} for related work.

\subsection*{Acknowledgements}
This paper arose from the author's response to a question of user \emph{rainmaker} on MathOverflow~\cite{dec-dioph}. It was originally
accepted for publication in~\cite{lib-albert}, but was omitted because of a technical mistake.

The research leading to these
results has received funding from the
European Research Council under the European Union's Seventh Framework
Programme (FP7/2007--2013)~/ ERC grant agreement no.~339691. The Institute of Mathematics of the Czech Academy of
Sciences is supported by RVO: 67985840.

I would like to thank the anonymous reviewers for many useful suggestions.

\section{Robinson defeats Diophantus}

We need a convenient way to refer to the individual axioms of~$Q$, thus we can as well start by properly defining the
theory, even though we trust it is familiar to the reader.
\begin{Def}\label{def:q}
$Q$ is the theory in language~$\LQ$ with axioms
\begin{gather}
\tag{Q1}Sx\ne0,\\
\tag{Q2}Sx=Sy\to x=y,\\
\tag{Q3}x=0\lor\exists y\,Sy=x,\\
\tag{Q4}x+0=x,\\
\tag{Q5}x+Sy=S(x+y),\\
\tag{Q6}x\cdot0=0,\\
\tag{Q7}x\cdot Sy=x\cdot y+x.
\end{gather}
Let $t\eqs u$ denote that the terms $t$ and~$u$ are syntactically identical. We define unary numerals $\num n\eqs S^n0$, and
binary numerals
\begin{align*}
\bnum 0&\eqs0,\\
\bnum{2n}&\eqs\num2\cdot\bnum n,\qquad n>0,\\
\bnum{2n+1}&\eqs S(\bnum{2n})
\end{align*}
for all natural numbers $n\in\N$. While unary numerals are easier to manipulate using
axioms of the theory, we will need the much shorter binary numerals when discussing algorithmic complexity. Of course,
$Q$ proves $\bnum n=\num n$, and we will use both interchangeably in contexts where the distinction does not matter.
\end{Def}

\subsection{Black-hole models}
As our starting point (already alluded to in the introduction), we can drastically reduce the complexity of the
Diophantine satisfiability problem for~$Q$ using black-hole models:
\begin{Lem}\label{lem:ninfty}
$\dio Q$ is polynomial-time reducible to $Q$-solvability of equations of the form $t=\bnum n$ with $n\in\N$.
\end{Lem}
\begin{Pf}
Consider the model $\N^\infty=\N\cup\{\infty\}$, where
$S\infty=\infty+x=x+\infty=x\cdot\infty=\infty$ for all $x\in\N^\infty$, $\infty\cdot0=0$, and $\infty\cdot x=\infty$
for $x\ne0$. It is readily seen that $\N^\infty\model Q$.

When written in binary, the lengths of $n+m$ and $n\cdot m$ are
bounded by the sum of lengths of $n$ and~$m$. It follows by induction on the complexity of~$t$ that given a term~$t$
and $\vec a\in\N^\infty$, the length of the value of $t(\vec a)$ in~$\N^\infty$ is polynomial in the lengths of $t$
and~$\vec a$, and we can compute $t(\vec a)$ in polynomial time.

Crucially, the operations in~$\N^\infty$ are defined so that they give a finite value only when forced so by the axioms
of~$Q$, hence we can show by induction on the complexity of~$t$ that
\[t(\vec\infty)=n\in\N\implies Q\vdash t(\vec x)=\bnum n.\]
For example, let $t\eqs u\cdot v$. Then $t(\vec\infty)\in\N$ only if both $u(\vec\infty),v(\vec\infty)\in\N$, or if
$v(\vec\infty)=0$. In the former case, the induction hypothesis gives
\[Q\vdash u(\vec x)=\bnum k,\qquad Q\vdash v(\vec x)=\bnum l\]
for some $k,l\in\N$, thus $Q\vdash t(\vec x)=\bnum n$ with $n=kl$. In the latter case, $Q\vdash v(\vec x)=0$ by the
induction hypothesis, hence $Q\vdash t(\vec x)=0$.

Thus, here is the promised reduction: given an equation $t_0=t_1$, if $t_0(\vec\infty)=\infty=t_1(\vec\infty)$, we have a
witness that $t_0=t_1$ is satisfiable, solving the problem outright;
otherwise, at least one of the terms $t_i$ is provably equal to a numeral~$\bnum n$, which we can compute in polynomial
time. The output of the reduction is (say) ``$0=0$'' in the former case, and ``$t_{1-i}=\bnum n$'' in the latter case.
\end{Pf}

This is not yet the end of the story; we can further reduce the problem by unwinding the terms from top. For example,
axioms Q1 and~Q2 imply that an equation $St=\num n$ is $Q$-satisfiable if and only if $n$ is nonzero, and
$t=\num{n-1}$ is satisfiable. Something to a similar effect also holds for the other function symbols, so let us see
where it gets us.
\begin{Def}
Let $Q_\forall$ denote the theory axiomatized by Q1, Q2, Q4--Q7, and
\begin{align}
\tag{Q8${}_n$}x+y=\num n&\to\LOR_{m\le n}(y=\num m),\\
\tag{Q9${}_n$}x\cdot y=\num n&\to x=0\lor\LOR_{m\le n}(y=\num m)
\end{align}
for $n\in\N$.
\end{Def}
\pagebreak[2]
\begin{Lem}\label{lem:univeasy}
\ \begin{enumerate}
\item\label{item:7}
$Q\vdash Q_\forall$.
\item\label{item:8}
Let $n\in\N$. Then $Q_\forall$ proves
\begin{align*}
x+y=\num n&\to\LOR_{k+m=n}(x=\num k\land y=\num m),\\
x\cdot y=\num n&\to x=0\lor\LOR_{km=n}(x=\num k\land y=\num m),\qquad n>0.
\end{align*}
\end{enumerate}
\end{Lem}
\begin{Pf}
\ref{item:7}: We prove Q8${}_n$ by induction on~$n$. Reason in~$Q$, and assume $x+y=\num n$. If $y=0$, we are done.
Otherwise $y=Sz$ for some~$z$ by Q3, hence $S(x+z)=\num n$ by Q5. This is only possible if $n>0$ due to Q1, and we
have $x+z=\num{n-1}$ by Q2, thus $z=0\lor\dots\lor z=\num{n-1}$ by the induction hypothesis, and consequently
$y=\num1\lor\dots\lor y=\num n$. (Alternatively, notice that under the traditional definition of $u\le v$ as $\exists
w\,(v=w+u)$, Q8${}_n$ may be read as the bounded sentence $\forall y\le\num n\,\LOR_{m\le n}y=\num m$,
hence its provability follows from the $\Sigma_1$-completeness of~$Q$.)

The proof of Q9${}_n$ is similar. Assuming $xy=\num n$, we are done if $y=0$, hence we can assume $y=Sz$. Then
$xz+x=\num n$ (Q7), thus $x=\num k$ for some $k=0,\dots,n$ by Q8${}_n$. If $k=0$, we are done. Otherwise 
$xz=\num{n-k}$ (Q4, Q5, Q2), where $n-k<n$, hence we can use the induction hypothesis to conclude $z=0\lor\dots\lor
z=\num{n-k}$. This implies $y=\num1\lor\dots\lor y=\num{n-k+1}$, where $n-k+1\le n$.

\ref{item:8}: If $x+y=\num n$, we have $y=\num m$ for some $m\le n$ by Q8${}_n$. Then $x+y=S^mx$ by Q4 and~Q5, hence
$x=\num{n-m}$ by~Q2.

Let $n\ne0$, and reason in~$Q_\forall$ again. Assume $xy=\num n$. We have $\num n\ne0$ by~Q1, hence $y\ne0$ by~Q6. Using
Q9${}_n$, either $x=0$, or $y=\num m$ for some $m=1,\dots,n$. In the latter case, $\num n=x\cdot\num{m-1}+x$ by~Q7,
hence $x=\num k$ for some $k=0,\dots,n$ by Q8${}_n$. Then $\num n=xy=\num{km}$ using Q4--7, hence $km=n$ by Q1, Q2.
\end{Pf}
\begin{Prop}\label{prop:q+}
Let $Q^+$ denote $Q$ extended by the axiom $0\cdot x=0$. Then $\dio{Q^+}$ is decidable.
\end{Prop}
\begin{Pf}
The proof of Lemma~\ref{lem:ninfty} works for~$Q^+$, too, with $\N^\infty$ modified so that $0\cdot\infty=0$. We
describe below a recursive procedure $\sol(E)$ that checks whether a finite set~$E$ of equations of the form $t=\bnum
n$ is $Q^+$-satisfiable.

Let $t=\bnum n$ be the first equation in~$E$ such that $t$ is not a variable, and $E'=E\bez\{t=\bnum n\}$:
\begin{enumerate}
\item\label{item:1} If $t\eqs t_0\cdot t_1$ and $n\ne0$, call $\sol(E'\cup\{t_0=\bnum{n_0},t_1=\bnum{n_1}\})$ for
every $n_0,n_1$ such that $n_0n_1=n$. Accept if any of the recursive calls accepted, otherwise reject.
\item\label{item:2} If $t\eqs t_0\cdot t_1$ and $n=0$, call $\sol(E'\cup\{t_0=0\})$ and $\sol(E'\cup\{t_1=0\})$. Accept if any of the
recursive calls accepted, otherwise reject.
\item\label{item:3} If $t$ is of the form $t_0+t_1$, $St_0$, or~$0$, proceed similarly.
\item\label{item:4} If the left-hand sides of all equations in~$E$ are variables, reject if $E$ contains a pair of equations with the
same left-hand sides and different right-hand sides, otherwise accept.
\end{enumerate}

Each recursive call strictly decreases the total number of symbols on the left-hand sides, hence the algorithm
terminates, and Lemma~\ref{lem:univeasy} and the extra axiom guarantee its correctness.

We note that $\sol(E)$ as presented is an exponential-time algorithm, but we can transform it into a nondeterministic
polynomial-time algorithm by making only one, nondeterministically chosen, recursive call at each step. Thus,
$\dio{Q^+}\in\np$.
\end{Pf}

If we try to use $\sol(E)$ for~$Q$, we run into trouble: while $Q$ proves $xy=0\to x=0\lor y=0$
by Lemma~\ref{lem:univeasy}, it does not prove the converse implication, hence the solvability of $E'\cup\{t_0=0\}$ does
not imply the solvability of $E'\cup\{t_0\cdot t_1=0\}$ in step~\ref{item:2}. Likewise, step~\ref{item:1} is incorrect,
because $E'\cup\{t_0\cdot t_1=\num n\}$ with $n\ne0$ may be satisfied in such a way that $t_0=0$.

However, the other reductions remain valid, and this still proves useful: a variant of $\sol(E)$ shows that $\dio Q$
reduces to $Q$-solvability of systems of equations of the form
\begin{equation}\label{eq:1}
\left\{\begin{aligned}
0\cdot t_1(\vec x)&=\num{n_1},\\
\cdots\span\\
0\cdot t_k(\vec x)&=\num{n_k}.
\end{aligned}\right.
\end{equation}
Diophantine systems of this type have not yet received the attention they deserve, so we are on our own. Their
$Q$-satisfiability turns out to be an unexpectedly circuitous problem: on the one hand, we will see that nearly every
such equation is satisfiable by itself in a suitable model, on the other hand there are subtle
dependencies that make systems such as
\begin{align*}
0\cdot(x+\num2)&=\num 5\\
0\cdot(y+0\cdot x)&=\num 7\\
0\cdot Sy&=\num 4
\end{align*}
unsatisfiable\footnote{By Q4--7, $0\cdot(x+\num2)=0\cdot SSx=(0\cdot x+0)+0=0\cdot x$, thus the first equation implies
$0\cdot x=\num5$. Likewise, the third equation gives $0\cdot y=\num4$, while the second equation gives
$\num7=0\cdot(y+\num5)=0\cdot y$. It is worth noting that the second equation does not imply
anything about $0\cdot y$ on its own---we need to know that $0\cdot x$ is standard first.}.
One consequence is that we cannot make do with a one-size-fits-all model of~$Q$ like in Lemma~\ref{lem:ninfty}; we
will need a variety of countermodels for different systems. This will be our task in the next two subsections.

\subsection{Universal fragment of $Q$}
We intend to use term models of a kind as our supply of models to satisfy various equations, but this approach
is not very friendly to the predecessor axiom~Q3, so to make our lives easier, we first determine what structures can
be \emph{extended} to models of~$Q$ by adding predecessors (and other elements that are forced upon us). By general
model theoretic considerations, these are exactly the models of the \emph{universal fragment} of~$Q$, hence we can
reformulate the problem as a description of this universal fragment. Since we used very suggestive notation,
the answer should come as no surprise:
\begin{Prop}\label{lem:quniv}
$Q_\forall$ is the universal fragment of~$Q$. That is, every model of~$Q_\forall$ embeds in a model of~$Q$.
\end{Prop}
\begin{Pf}
Fix $M\model Q_\forall$. Identifying each $n\in\N$ with the corresponding numeral $\num n^M\in M$, we may assume that
$M$ includes the standard model~$\N$; in particular, $M\bez\N$ is the set of nonstandard elements of~$M$.

Let $A$ denote the set of nonzero elements of~$M$ without a predecessor. We will embed $M$ in
a structure with domain
\[N=M\cup\{\infty\}\cup\{\p{a,k},\p{a,k,x}:a\in A,k\in\N^{>0},x\in M\bez\N\},\]
where $\p{a,k}$ should be thought of as $a-k$, and $\p{a,k,x}$ as $x\cdot(a-k)$. We will define the interpretations of
the $\LQ$-function symbols in~$N$, and verify that $N\model Q$ along the way. All the function symbols are understood
to retain their interpretations from~$M$ on elements of~$M$, thus we will not indicate such cases explicitly.

\emph{Successor:} we put
\begin{align*}
S^N\infty&=\infty,\\
S^N\p{a,k,x}&=\p{a,k,x},\\
S^N\p{a,k}&=\begin{cases}\p{a,k-1}&k>1,\\a&k=1.\end{cases}
\end{align*}
We can see immediately that this makes $N$ a model of Q1--Q3. This also means we can unambiguously refer to $S^nx$
for $n\in\Z$ and $x\in N\bez\N$.

\emph{Addition:} we put
\begin{align*}
\infty+^Ny&=\infty,\\
\p{a,k}+^Ny&=\begin{cases}S^{n-k}a&y=n\in\N,\\\infty&\text{otherwise.}\end{cases}
\intertext{For $x\in M$, we define}
x+^N\infty&=\infty,\\
x+^N\p{a,k,y}&=\infty,\\
x+^N\p{a,k}&=S^{-k}(x+^Ma).
\intertext{Note that the last item makes sense: since $a\notin\N$, also $x+^Ma\notin\N$ by~Q8. Finally, we put}
\p{a,k,x}+^Nn&=\p{a,k,x},\qquad n\in\N,\\
\p{a,k,x}+^NS^nx&=\begin{cases}\p{a,k-1,x}&k>1,\\S^n(x\cdot^Ma)&k=1,\end{cases}\qquad n\in\Z,\\
\p{a,k,x}+^Ny&=\infty\qquad\text{for other $y$.}
\end{align*}
Again, $x\cdot^M a\notin\N$ by Q9 as $x,a\notin\N$.

It is straightforward to check that $N$ validates Q4 and~Q5.

\emph{Multiplication:} for $x\in N\bez M$, we put
\begin{align*}
x\cdot^Nn&=(\cdots(0+^N\underbrace{x)+^N\dots+^Nx)+^Nx}_n,\qquad n\in\N,\\
x\cdot^Ny&=\infty,\qquad y\notin\N.
\intertext{For $x\in M$, we define}
x\cdot^N\infty&=\infty,\\
x\cdot^N\p{a,k,y}&=\infty,\\
x\cdot^N\p{a,k}&=\begin{cases}S^{-kn}(n\cdot^Ma)&x=n\in\N,\\\p{a,k,x}&x\notin\N.\end{cases}
\end{align*}
As above, $S^{-kn}(n\cdot^Ma)$ exists: either $n=0$ and the $S^{-kn}$ does nothing, or $n>0$, in which case
$n\cdot^Ma\notin\N$ by Q9.

Again, it is straightforward to check Q6 and~Q7.
\end{Pf}

\subsection{Reduced terms}
We now come to the crucial part of our construction: we establish that a system~\eqref{eq:1} is $Q$-satisfiable if the
terms~$t_i$ obey certain conditions that guarantee they do not interact with each other.

\begin{Def}\label{def:norm}
An $\LQ$-term is \emph{normal} if it contains no subterm of the form $t+0$, $t+Su$, $t\cdot0$, or $t\cdot Su$.
A normal term is \emph{irreducible} if it does not have the form $0$ or~$St$. In other words, normal and irreducible
terms are generated by the following grammar:
\begin{align*}
I&::=x_i\mid (N+I)\mid (N\cdot I)\\
N&::=I\mid0\mid SN
\end{align*}
If $T$ is a set of terms, a normal term is \emph{$T$-reduced} if it contains no subterm of the form $0\cdot t$
for $t\in T$.

Each normal term can be uniquely written in the form $S^n0$ or $S^nt$, where $n\in\N$, and $t$ is irreducible. A
subterm of a normal ($T$-reduced) term is again normal ($T$-reduced, resp.).
\end{Def}

\begin{Lem}\label{lem:normal}
Let $T=\{t_i:i<k\}$ be a finite sequence of distinct irreducible terms such that $0\cdot t_i$ is not a subterm of~$t_j$ for any
$i,j<k$ (i.e., the terms $t_i$ are $T$-reduced), and $\{n_i:i<k\}\sset\N$.
\begin{enumerate}
\item\label{item:9}
The set of equations $\{0\cdot t_i=\num{n_i}:i<k\}$ is $Q$-satisfiable.
\item\label{item:10}
Given $T,\vec n$, and a term~$t$, we can compute a $T$-reduced term $\red t$ such that
\[Q\vdash\ET_{i<k}0\cdot t_i(\vec x)=\num{n_i}\to t(\vec x)=\red t(\vec x).\]
\end{enumerate}
\end{Lem}
\begin{Pf}
\ref{item:9}: We define a model~$M$ whose domain consists of all $T$-reduced terms, and operations as follows. We put
$0^M=0$, and $S^Mt=St$. If $t\in M$, $n\in\N$, and $u\in M$ is irreducible,
\begin{align*}
t+^MS^n0&=S^nt,\\
t+^MS^nu&=S^n(t+u).
\intertext{If $t,u\in M$ are irreducible, and $n,m\in\N$, we put}
S^nt\cdot^MS^m0&=\underbrace{S^n(S^n(\dots(S^n}_m(0+\underbrace{t)\dots)+t)+t}_m),\\
S^nt\cdot^MS^mu&=\underbrace{S^n(S^n(\dots(S^n}_m(S^nt\cdot u+\underbrace{t)\dots)+t)+t}_m),\\
S^n0\cdot^MS^m0&=S^{nm}0,\\
S^n0\cdot^MS^mu&=S^{nm}(S^n0\cdot u),\qquad n>0,\\
0\cdot^MS^mu&=\begin{cases}S^{n_i}0&u=t_i,\\0\cdot u&u\notin T.\end{cases}
\end{align*}

It is readily checked that the operations are well-defined (i.e., the terms given above as their values are
$T$-reduced), and that $M\model Q_\forall$. Let $v$ be the valuation in~$M$ which assigns each variable~$x_i$ to the
corresponding element $x_i\in M$. Then $v(t)=t$ for every $T$-reduced term~$t$, hence $v$ satisfies in~$M$ the equations
$0\cdot t_i=\num{n_i}$. By Proposition~\ref{lem:quniv}, we can embed~$M$ into a model of~$Q$.

\ref{item:10}: Since the operations in~$M$ are computable, we can compute the value $v(t)\in M$ by induction on the
complexity of~$t$. This value is a $T$-reduced term, so we can define $\red t=v(t)$. Then we show that $Q$ proves the
required implication
\begin{equation}\label{eq:2}
\ET_{i<k}0\cdot t_i(\vec x)=\num{n_i}\to t(\vec x)=\red t(\vec x)
\end{equation}
by induction on the complexity of~$t$. For the induction steps, we observe the operations in~$M$ are defined so that if
$t+^Ms=u$, then $Q\vdash t+s=u$, and likewise for~$\cdot^M$ with the exception of the clause $0\cdot^MS^mt_i=S^{n_i}0$,
which is handled by the premise of~\eqref{eq:2}.
\end{Pf}

The reader might have realized that what just happened was term rewriting in thinly veiled disguise. Even though we
will not need this point of view for our application, we make the digression to spell this connection out because of sheer curiosity.
\begin{Def}\label{def:rs}
Let $R_Q$ denote the rewriting system for $L_Q$-terms generated by the rules
\begin{equation}\label{eq:4}
\left\{\quad\begin{aligned}
t+0&\lead t,\\
t+Su&\lead S(t+u),\\
t\cdot0&\lead0,\\
t\cdot Su&\lead t\cdot u+t.
\end{aligned}\right.
\end{equation}
More generally, if $\{t_i:i<k\}$ is a sequence of terms satisfying
the conditions of Lemma~\ref{lem:normal}, and $\{n_i:i<k\}\sset\N$, let $R_{\vec t,\vec n}$ denote the rewriting system
extending $R_Q$ with the rules 
\begin{equation}\label{eq:3}
0\cdot t_i\lead \num{n_i},\qquad i<k
\end{equation}
(these rules are not supposed to allow substitution for variables inside~$t_i$).
\end{Def}

Notice that a term is normal in the sense of Definition~\ref{def:norm} iff it is a normal form with respect to~$R_Q$,
and it is $T$-reduced (with $T=\{t_i:i<k\}$) iff it is a normal form with respect to $R_{\vec t,\vec n}$ for an
arbitrary choice of~$\vec n$.
\begin{Prop}\label{lem:trs}
For any $\vec t$ and~$\vec n$ as in the definition, the rewriting system $T_{\vec t,\vec n}$ is strongly normalizing
and confluent. (That is, every term has a unique normal form, and every sequence of reductions will eventually reach
it.)
\end{Prop}
\begin{Pf}
Put $c=2+\max_{i<k}n_i$, and define a ``norm'' function on terms by
\begin{align*}
\nrm{x_i}=\nrm0&=c,\\
\nrm{St}&=\nrm t+3,\\
\nrm{t+u}&=\nrm t+2\nrm u,\\
\nrm{t\cdot u}&=\nrm t\cdot\nrm u.
\end{align*}
Notice that $\nrm t\ge c$ for any term~$t$, and the norm is strictly monotone in the sense that $\nrm u<\nrm v$ implies
$\nrm{t(u)}<\nrm{t(v)}$. Using this, we can check easily that all $R_{\vec t,\vec n}$-reduction steps strictly decrease
the norm, thus there is no infinite sequence of reductions: in particular, we have
\[\nrm{\num{n_i}}=3n_i+c\le 4c-6<c^2\le\nrm{0\cdot t_i}.\]
This shows strong normalization of~$R_{\vec t,\vec n}$.

By Newman's lemma, confluence is implied by local confluence: that is, it suffices to show that if $s\lead v_0$ and
$s\lead v_1$, then $v_0\leadt w$ and $v_1\leadt w$ for some term~$w$, where $\lead$ denotes one-step reduction, and
$\leadt$ its reflexive transitive closure.

The local confluence property obviously holds if the two reductions $s\lead v_i$ are identical, or if they operate on
disjoint terms. It also holds if $s\lead v_i$ is one of the $R_Q$-reductions as given in~\eqref{eq:4}, and
$s\lead v_{1-i}$ operates inside one of the terms $t,u$ on the left-hand side of~\eqref{eq:4}: we can
instead perform the reduction on their copies on the right-hand side.

This in fact covers all possibilities: the only redexes properly included inside the left-hand side of any $R_Q$-rule
in~\eqref{eq:4} are inside $t$ or~$u$, as there are no rules reducing $0$ or $Su$; there are no
redexes properly included inside $0\cdot t_i$ by the assumption that $t_i$ is $T$-reduced; and each redex can be
reduced only in one way---the only possible clashes could be between \eqref{eq:3} and the $R_Q$-rules for multiplication, but
these are prevented as $t_i$ is assumed not to be of the form $0$ or~$Su$.
\end{Pf}

Proposition~\ref{lem:trs} provides an alternative proof for most of Lemma~\ref{lem:normal}: first, the $T$-reduced term~$\red
t$ in \ref{lem:normal}~\ref{item:10} is just the $R_{\vec t,\vec n}$-normal form for~$t$. Second, we can use confluence (Church--Rosser
property) to construct the model $M$ for~\ref{item:9} as the model of $R_{\vec t,\vec n}$-normal terms, or
equivalently, as the quotient of the free term model by the equivalence relation induced by reduction. It is
automatically a model of axioms Q4--7 embodied in the reduction rules, and it is easily seen to satisfy Q1 and~Q2
because there are no rules with redex~$Su$. It would still take a little work to establish the validity of Q8 and~Q9.

\subsection{Witnessing satisfiability}
To complete our analysis of~$\dio Q$, we will now show that a general equation $t=\num n$ can only be $Q$-satisfied if it is
implied by a (suitably bounded) system of the form~\eqref{eq:1} that respects the assumptions of Lemma~\ref{lem:normal}.

\pagebreak[2]
\begin{Def}\label{def:wit}
For any term~$u$, let $\red u$ denote its $\nul$-reduced form as given by Lemma~\ref{lem:normal}.

A \emph{labelling} of a term~$t$ is a partial map~$\ell$ from subterms of~$t$ to~$\N$. If $\ell$ is a
labelling of~$t$, and $u$ a subterm of~$t$ (written henceforth as $u\sset t$), let $u_\ell$ be the term obtained from $u$ by replacing all maximal proper labelled
subterms of~$u$ by numerals for their labels.

A \emph{witness} for $t=\num n$ is a labelling $\ell$ of~$t$ by numbers $k\le n$ such that:
\begin{enumerate}
\item\label{item:14} $\ell(t)=n$.
\item\label{item:15}
If $u,v\sset t$ are such that $\red{u_\ell}\eqs\red{v_\ell}$, then $\ell(u)=\ell(v)$ (meaning both are undefined, or
both are defined and equal).
\item\label{item:16}
If $u\in\dom(\ell)$, and $\red{u_\ell}\eqs\num k$ for some~$k\in\N$, then $\ell(u)=k$.
\item\label{item:17}
If $u\in\dom(\ell)$, then all immediate subterms of~$u$ are labelled, unless $u\eqs v\cdot w$, and $v$ or~$w$ is
labelled~$0$.
\end{enumerate}
\end{Def}
Note that \ref{item:15} implies that occurrences of the same subterm either all have the same label, or are all
unlabelled. We also remark that in~\ref{item:16}, $k\le n$ is not a premise, but part of the conclusion.
\begin{Exm}
Table~\ref{tab:label} shows a labelling~$\ell$ of the term $t\eqs x\cdot y+x\cdot SSSy$ that is a witness for satisfiability
of the equation $t=\num8$. For convenience, the table also lists for each term $u\sset t$ its set of maximal proper
labelled subterms, as well as $u_\ell$ and $\red{u_\ell}$, which makes it easy to check that conditions
\ref{item:14}--\ref{item:17} hold. In particular, for~\ref{item:15}, the two terms $u$ with $\red{u_\ell}\eqs 0\cdot y$ have the same
label $\ell(u)=4$; for~\ref{item:16}, the only applicable case is $\red{t_\ell}=\num8$, which agrees with $\ell(t)=8$.
We invite the reader to verify that $\ell$ is in fact the \emph{only} possible witness for $t=\num8$.

For this example, the set~$E$ considered below in the proof of Lemma~\ref{lem:wit} consists of the single equation
$0\cdot y=\num4$.
\end{Exm}
\begin{table}
\[\begin{array}{ccccc}
u&\ell(u)&\text{m.p.l.s.}&u_\ell&\red{u_\ell}\\
\hline
x\cdot y+x\cdot SSSy&8&x\cdot y,x\cdot SSSy&\num4+\num4&\num8\\
x\cdot y&4&x&0\cdot y&0\cdot y\\
x\cdot SSSy&4&x&0\cdot SSSy&0\cdot y\\
x&0&\text{--}&x&x\\
S^iy\ (i=0,\dots,3)&\text{--}&\text{--}&S^iy&S^iy
\end{array}\]
\caption{Witness for $x\cdot y+x\cdot SSSy=\num8$}
\label{tab:label}
\end{table}
\pagebreak[2]
\begin{Lem}\label{lem:wit}
An equation $t=\num n$ is $Q$-satisfiable if and only if it has a witness.
\end{Lem}
\begin{Pf}
Left-to-right: let $M\model Q$ and $\vec a\in M$ be such that $t^M(\vec a)=n$. Define a labelling of~$t$ by putting
$\ell(u)=u^M(\vec a)$ if $u^M(\vec a)\in\{0,\dots,n\}$, and $\ell(u)$ is undefined otherwise. Since a term equals
its $\nul$-reduction provably in~$Q$, we have $\red{u_\ell}^M(\vec a)=u^M(\vec a)$ for any $u\sset t$. It follows
easily that $\ell$ is a witness for $t=\num n$, using Lemma~\ref{lem:univeasy} for condition~\ref{item:17}.

Right-to-left: let $E$ denote the set of equations
\[\red{u_\ell}=\num k\]
where $u\eqs v\cdot w\sset t$, $\ell(u)=k$, $\ell(v)=0$, and $\red{u_\ell}\neqs0$ (which implies $w\notin\dom(\ell)$). Note that
$\red{u_\ell}$ then must be of the form $0\cdot u^-$, where $u^-$ is an irreducible term, and
$\red{w_\ell}\eqs S^mu^-$ for some~$m$.
\pagebreak[2]
\begin{Cl}
If $E$ is satisfiable, then $t=\num n$ is satisfiable.
\end{Cl}
\begin{Pf*}
Fix a model $M\model Q$ and $\vec a\in M$ that satisfies~$E$. Note that $E$ only contains unlabelled
variables\footnote{Variables in $\red{u_\ell}$ come from variables in~$u$. However, a
labelled variable in~$u$ is a maximal proper labelled subterm of~$u$, or is included in such a maximal
subterm; consequently, it disappears in $u_\ell$ (and $\red{u_\ell}$) by virtue of being replaced with a constant term
(a numeral).}; if
$x_i$ is labelled, we make sure that $a_i=\ell(x_i)$ (this is independent of the choice of an occurrence of $x_i$
in~$t$ by condition~\ref{item:15}). We claim that
\[u\in\dom(\ell)\implies u^M(\vec a)=\ell(u),\]
which gives $t^M(\vec a)=n$ by condition~\ref{item:14}. We prove this by induction on the complexity of~$u$.

The statement holds for variables, and condition~\ref{item:16} implies it holds for $u\eqs 0$.

If $u\in\dom(\ell)$ is of the form $Sv$ or $v+w$, then $v,w\in\dom(\ell)$ by~\ref{item:17}, and $\ell(u)$ equals $\ell(v)+1$ or
$\ell(v)+\ell(w)$ (resp.) by~\ref{item:16}, thus $u^M(\vec a)=\ell(u)$ by the induction hypothesis for $v$ and~$w$.

The same argument applies if $u\eqs v\cdot w$, and both $v,w\in\dom(\ell)$, or $\ell(w)=0$. Assume $\ell(v)=0$ and
$w\notin\dom(\ell)$. Using the induction hypothesis for subterms of~$u$, and the soundness of reduction, we
have $u^M(\vec a)=u_\ell^M(\vec a)=\red{u_\ell}^M(\vec a)$. If $\red{u_\ell}=0$, this means $u^M(\vec a)=0=\ell(u)$
by~\ref{item:16}; otherwise the equation $\red{u_\ell}=\ell(u)$ is in~$E$, hence it is satisfied by~$\vec a$.
\end{Pf*}
Condition~\ref{item:15} ensures that $E$ does not contain two equations with the same left-hand side. Moreover, for any
\begin{align*}
0\cdot u^-_0&=\num{k_0}\\
0\cdot u^-_1&=\num{k_1}
\end{align*}
in~$E$, $0\cdot u^-_0$ is not a subterm of $u^-_1$: writing $u_1\eqs v_1\cdot w_1$, inspection of the definition of reduction
shows that this could only happen if $(w_1)_\ell$ contained a nonconstant subterm $s$ such that $\red s\eqs 0\cdot u^-_0$. But
then $s\eqs r_\ell$ for some $r\sset w_1$ such that $r\notin\dom(\ell)$, whereas we should have $\ell(r)=k_0$
by~\ref{item:15}, a contradiction. Thus, $E$ is satisfiable by Lemma~\ref{lem:normal}.
\end{Pf}

\begin{Thm}\label{thm:diophq}
$\dio Q$ is decidable.
\end{Thm}
\begin{Pf}
By Lemma~\ref{lem:ninfty}, $\dio Q$ reduces to $Q$-satisfiability of equations of the form $t=\num n$. These can be
checked by the criterion from Lemma~\ref{lem:wit}: a witness for $t=\num n$ has size bounded by a computable function
of $t$ and~$n$, and using the computability of~$\red u$, we can algorithmically recognize a witness when we see it.
\end{Pf}

\section{Computational complexity}

Our arguments thus far give an exponential-time algorithm for checking if a given Diophantine equation is $Q$-satisfiable. We can in fact determine the complexity of~$\dio Q$ precisely.
First, a general lower bound follows from a beautiful result of Manders and Adleman~\cite{man-adl} that there are very simple
$\np$-complete Diophantine problems.
\begin{Thm}[Manders and Adleman]\label{thm:ma}
The following problem is $\np$-complete: given $a,b\in\N$ in binary, determine whether $x^2+ay-b=0$ has a solution
in~$\N$.
\noproof\end{Thm}
(They state it with $ax^2+by-c$, but it is easy to show that the version here is equivalent.)

\pagebreak[2]
\begin{Cor}\label{cor:np-hard}
If $T$ is a consistent extension of~$Q_\forall$, then $\dio T$ is $\np$-hard.
\end{Cor}
\begin{Pf}
If $a>0$ (which we can assume without loss of generality), $x^2+ay-b=0$ is solvable iff the equation 
\begin{equation}\label{eq:7}
x\cdot x+\bnum a\cdot y=\bnum b
\end{equation}
is in~$\dio T$: on the one hand, a solution in~$\N$ yields a
solution in any model of~$T$. On the other hand, \eqref{eq:7} implies in~$Q_\forall$ that $x\cdot x$ and $\bnum a\cdot
y$ are standard and bounded by~$b$ using~Q8, hence $y$ is standard by~Q9. Also by~Q9, $x=0$, or $x=0,\dots,\bnum
b$; either way, $x$ is  standard. Thus, if \eqref{eq:7} is solvable in any model of~$T\Sset Q_\forall$, it is
solvable in~$\N$.
\end{Pf}

We will show that $\dio Q$ is as easy as possible, i.e., $\np$-complete. Now, the witnesses for satisfiability from
Definition~\ref{def:wit} are polynomial-size objects (if we write all numbers in binary), but it is not immediately
clear they can be recognized in polynomial time. In particular, the conditions demand us to test
$\red{u_\ell}\eqs\red{v_\ell}$ for subterms $u,v\sset t$, which na\"\i vely takes exponential time as the $\red
t$~reduction from Lemma~\ref{lem:normal} can exponentially blow up sizes of terms (e.g., it unwinds a binary
numeral term to the corresponding unary numeral). Fortunately, the offending overlarge pieces have a
very boring, repetitive structure, hence we can overcome this obstacle by devising a succinct representation
of terms such that on the one hand, the reduction of a given term has a polynomial-size representation, and on
the one hand, we can efficiently test whether two representations describe the same term.

The representations we use below (called \emph{descriptors}) have the syntactic form of terms over the language~$\LQ$
augmented with extra function symbols $S_n(x)$, $A_{n,m}(x)$, and $B_{n,m}(x,y)$, where $n,m$ are integer indices
written in binary. Their exact meaning is explained below, however, the intention is that they facilitate implementation of
the operations (especially multiplication) introduced in the proof of Lemma~\ref{lem:normal}.
\begin{Def}\label{def:desc}
We define a set of expressions called \emph{(term) descriptors}, and for each descriptor~$t$ a term $\des(t)$ which it
denotes, as follows.
\begin{itemize}
\item The constant~$0$ and variables~$x_i$ are descriptors denoting themselves.
\item If $t,u$ are descriptors, then $t+u$ and $t\cdot u$ are descriptors, and $\des(t+u)=\des(t)+\des(u)$,
$\des(t\cdot u)=\des(t)\cdot\des(u)$.
\item If $t$ is a descriptor, and $n\ge1$ is written in binary, then $S_n(t)$ is a descriptor, and
$\des(S_n(t))=S^n(\des(t))$.
\item If $u$ is a descriptor, and $n\ge0$, $m\ge2$ are written in binary, then $A_{n,m}(t)$ is a descriptor, and
\[\des(A_{n,m}(u))=\underbrace{S^n(\dots(S^n}_{m-1}(0+\underbrace{\des(u))\dots)+\des(u))+\des(u)}_m.\]
\item If $t,u$ are descriptors, and $n\ge0$, $m\ge1$ are written in binary, then $B_{n,m}(t,u)$ is a descriptor, and
\[\des(B_{n,m}(t,u))=\underbrace{S^n(\dots(S^n}_{m-1}(S^n(\des(u))\cdot\des(t)+\underbrace{\des(u))\dots)+\des(u))+\des(u)}_m.\]
\end{itemize}
A descriptor is \emph{minimal} if it contains no subdescriptors of the form
\begin{align*}
&S_n(S_m(t)),&&n,m\ge1,\\
&S_n(0+u)+u,&&n\ge0,\\
&S_n(u)\cdot t+u,&&n\ge0,\\
&S_n(A_{n,m}(u))+u,&&n\ge0,\:m\ge2,\\
&S_n(B_{n,m}(t,u))+u,&&n\ge0,\:m\ge1,
\end{align*}
where $S_0(t)$ is understood as~$t$.
\end{Def}

Notice that the definitions of $\des(A_{n,m}(u))$ and~$\des(B_{n,m}(t,u))$ are short of an outer~$S^n$ as compared to
the relevant clauses in Lemma~\ref{lem:normal}. The reason for this choice is that in the inductive construction
of~$\red t$, we need to be able to peel off easily the outer stack of $S$'s from the terms we got from the inductive
hypothesis in order to proceed.
\begin{Lem}\label{lem:desc}
\
\begin{enumerate}
\item\label{item:20}
Given a descriptor~$t$, we can compute in polynomial time a minimal descriptor~$t'$ such that
$\des(t)\eqs\des(t')$.
\item\label{item:19}
Given a term~$t$, we can compute in polynomial time a descriptor~$t'$ such that $\des(t')\eqs\red t$.
\item\label{item:18}
If $t_0,t_1$ are minimal descriptors such that $\des(t_0)\eqs\des(t_1)$, then $t_0\eqs t_1$.
\item\label{item:21}
Given descriptors $t$ and~$u$, we can test in polynomial time whether $\des(t)\eqs\des(u)$.
\end{enumerate}
\end{Lem}
\begin{Pf}
\ref{item:20}: We minimize the descriptor by applying the following rules to its subdescriptors in arbitrary
order:
\begin{align*}
S_n(S_m(t))&\lead S_{n+m}(t),\\
S_n(0+u)+u&\lead A_{n,2}(u),\\
S_n(u)\cdot t+u&\lead B_{n,1}(t,u),\\
S_n(A_{n,m}(u))+u&\lead A_{n,m+1}(u),\\
S_n(B_{n,m}(t,u))+u&\lead B_{n,m+1}(t,u),
\end{align*}
where $n,m$ are as appropriate for each case according to Definition~\ref{def:desc}. Each rule strictly decreases the
number of function symbols in the descriptor, hence the procedure stops after polynomially many steps, and it clearly
produces a minimal descriptor. Also, the maximal length (in binary) of numerical indices increases by at most~$1$ in
each step, hence all descriptors produced during the process have polynomial size, and the algorithm runs in polynomial
time.

\ref{item:19}: By a straightforward bottom-up approach mimicking the definition in Lemma~\ref{lem:normal}, we compute for each subterm $u\sset t$ a descriptor~$u'$
such that $\des(u')=\red u$, and $u'$ has the form $S_n(u'')$ where $u''$ is~$0$ or denotes an irreducible term. As in
Lemma~\ref{lem:ninfty}, all numerical indices appearing during the computation have bit-length bounded by the size
of~$t$. If $u\eqs u_0+u_1$ or $u\eqs u_0\cdot u_1$, then $u'$ can be expressed by at most one occurrence of each of
$u_0''$, $u_1''$, and a bounded number of other symbols (by employing the $A_{n,m}$ and~$B_{n,m}$ functions).
It follows easily that all descriptors constructed during the computation have polynomial size, and the computation
works in polynomial time.

\ref{item:18}: By induction on the complexity of $\des(t_0),\des(t_1)$. If $t_0\eqs u_0\cdot v_0$, then $t_1$ must be
of the form $u_1\cdot v_1$, as other descriptors denote terms whose topmost symbols are different from~$\cdot$. Then
$\des(u_0)\eqs\des(u_1)$ and $\des(v_0)\eqs\des(v_1)$, hence $u_0\eqs u_1$ and $v_0\eqs v_1$ by the induction
hypothesis, hence $t_0\eqs t_1$. A similar argument applies when the topmost symbol of $t_0$ or~$t_1$ is a variable, $0$,
or~$S_n$ (in the last case, we use the fact that if $t_i\eqs S_n(u_i)$, then $u_i$ cannot have topmost symbol~$S_m$ by
minimality).

The remaining cases are when both $t_i$ are of the forms $u_i+v_i$, $A_{n_i,m_i}(v_i)$, or
$B_{n_i,m_i}(u_i,v_i)$, so that the topmost symbol of $\des(t_i)$ is~$+$. We have $\des(t_i)\eqs\des(w_i)+\des(v_i)$, where
\[w_i\eqs
\begin{cases}
u_i&t_i\eqs u_i+v_i,\\
S_{n_i}(A_{n_i,m_i-1}(v_i))&t_i\eqs A_{n_i,m_i}(v_i),m_i\ge3,\\
S_{n_i}(0+v_i)&t_i\eqs A_{n_i,2}(v_i),\\
S_{n_i}(B_{n_i,m_i-1}(u_i,v_i))&t_i\eqs B_{n_i,m_i}(u_i,v_i),m_i\ge2,\\
S_{n_i}(v_i)\cdot u_i&t_i\eqs B_{n_i,1}(u_i,v_i).
\end{cases}
\]
Here $n_i$ may be~$0$, in which case $S_{n_i}$ is void. The descriptor $w_i$ as given here is minimal,
except that in the last case, it might happen that $v_i$ itself starts with $S_{k_i}$ for some~$k_i$; in that case, we modify $w_i$ in
the obvious way. Since $\des(w_0)\eqs\des(w_1)$ and $\des(v_0)\eqs\des(v_1)$ are proper subterms of
$\des(t_0)\eqs\des(t_1)$, we may now apply the induction hypothesis, yielding $v_0\eqs v_1$, and $w_0\eqs w_1$. By
inspection, we see that for each of the five clauses of the definition of~$w_i$, we can read off the original
parameters ($n_i$, $m_i$, $u_i$; we already know $v_i$) from~$w_i$. Moreover, two distinct clauses cannot result in the
same~$w_i$: the only problematic case is the first clause, where we need to use the minimality of~$t_i$. Thus, all in
all, $w_i$ and~$v_i$ uniquely determine~$t_i$, hence we obtain $t_0\eqs t_1$.

\ref{item:21} follows from \ref{item:20} and~\ref{item:18}.
\end{Pf}

\pagebreak[2]
\begin{Thm}
$\dio Q$ is $\np$-complete.
\end{Thm}
\begin{Pf}
$\np$-hardness is Corollary~\ref{cor:np-hard}, hence in view of Lemmas \ref{lem:ninfty} and~\ref{lem:wit}, it suffices to show that we
can check the existence of a witness~$\ell$ for $t=\bnum n$ in~$\np$. It is immediate from the definition that $\ell$
has size polynomial in $\log n$ and in the length of~$t$ if we write labels in binary, so it remains to verify conditions
\ref{item:14}--\ref{item:17} in polynomial time.

Conditions \ref{item:14} and~\ref{item:17} are clearly polynomial-time. As for~\ref{item:15}, notice first that it
makes no difference whether we use unary or binary numerals in the construction of~$u_\ell$, as both end up the same
after applying~$\red\ $. Thus, in order to test $\red{u_\ell}\eqs\red{v_\ell}$ in polynomial time, we can compute
$u_\ell,v_\ell$ using binary numerals, compute descriptors denoting $\red{u_\ell},\red{v_\ell}$ using
Lemma \ref{lem:desc}~\ref{item:19}, and compare them using Lemma \ref{lem:desc}~\ref{item:21}.

Condition~\ref{item:16} is similar: given a term~$u$, we can compute a minimal descriptor for $\red{u_\ell}$ in polynomial time, and then check easily whether it has the form $S_k(0)$, and
if so, extract~$k$.
\end{Pf}

\section{Conclusion}

Unlike stronger theories of arithmetic, we have seen that satisfiability of Diophantine equations in models of~$Q$ can
be tested in~$\np$, hence undecidability only sets in for more complicated $\Sigma_1$~sentences. The proof also
revealed that Robinson's arithmetic can divide standard numbers by zero with ruthless efficiency (albeit in a lopsided
way).

Some related questions suggest themselves, such as how far can we push the argument? On the one hand, the
criterion in Lemma~\ref{lem:wit} does not use in any way that we are dealing with a single equation. Considering also that
the models constructed in Lemma~\ref{lem:normal} only equate terms with the same reduced form, we obtain easily the
following generalization:
\begin{Prop}\label{prop:gener}
$Q$-satisfiability of existential sentences, all of whose positively occurring atomic subformulas are of the form
$t=\bnum n$, is decidable, and $\np$-complete.
\noproof\end{Prop}
On the other hand, the reduction in Lemma~\ref{lem:ninfty} breaks down already for conjunctions of two equations, hence we
are led to
\begin{Prob}\label{prob:exist}
Is $Q$-satisfiability of existential sentences decidable?
\end{Prob}

A question in another vein is how much stronger can we make the theory while maintaining decidability. Observe that the
simple argument in Proposition~\ref{prop:q+} applies not just to~$Q^+$ itself, but also to all its extensions valid in the
variant $\N^\infty$~model used in the proof. This model is actually quite nice: a totally ordered commutative
semiring, one pesky axiom short of the theory~$\PA^-$!
\begin{Prob}\label{prob:pa-}
Is $\dio{\PA^-}$ decidable?
\end{Prob}
This problem appears to be essentially as hard as the decidability of $\dio\io$ mentioned in the introduction, cf.\
\cite{wil:iop,vdd:dor}.


\begin{thebibliography}{10}

\bibitem{doy-con}
Peter~G. Doyle and John~H. Conway, \emph{Division by three},
  {arXiv:math/0605779 [math.LO]}, 1994,
  \url{http://arxiv.org/abs/math/0605779}.

\bibitem{doy-qiu}
Peter~G. Doyle and Cecil Qiu, \emph{Division by four}, {arXiv:1504.01402
  [math.LO]}, 2015, \url{http://arxiv.org/abs/1504.01402}.

\bibitem{vdd:dor}
Lou van~den Dries, \emph{Which curves over {$\mathbf Z$} have points with
  coordinates in a discrete ordered ring?}, Transactions of the American
  Mathematical Society 264 (1981), no.~1, pp.~181--189.

\bibitem{lib-albert}
Jan van Eijck, Rosalie Iemhoff, and Joost~J. Joosten (eds.), \emph{Liber
  amicorum Alberti: A tribute to Albert Visser}, Tributes vol.~30,
  College Publications, London, 2016.

\bibitem{mrdpea}
Haim Gaifman and Constantinos Dimitracopoulos, \emph{Fragments of Pe\-ano's
  arithmetic and the {MRDP} theorem}, in: Logic and algorithmic, Monographie de
  L'En\-sei\-gne\-ment Ma\-th{\'e}\-ma\-tique no.~30, Universit{\'e} de
  Gen{\`e}ve, 1982, pp.~187--206.

\bibitem{kaye:dio}
Richard Kaye, \emph{Diophantine induction}, Annals of Pure and Applied Logic 46
  (1990), no.~1, pp.~1--40.

\bibitem{kaye:hilb}
\bysame, \emph{Hilbert's tenth problem for weak theories of arithmetic},
  Annals of Pure and Applied Logic 61 (1993), no.~1--2, pp.~63--73.

\bibitem{book}
Jan Kraj{\'\i}{\v c}ek, \emph{Bounded arithmetic, propositional logic, and
  complexity theory}, Encyclopedia of Mathematics and Its Applications vol.~60,
  Cambridge University Press, 1995.

\bibitem{man-adl}
Kenneth~L. Manders and Leonard~M. Adleman, \emph{$\mathit{NP}$-complete
  decision problems for binary quadratics}, Journal of Computer and System
  Sciences 16 (1978), no.~2, pp.~168--184.

\bibitem{otero}
Margarita Otero, \emph{On Diophantine equations solvable in models of open
  induction}, Journal of Symbolic Logic 55 (1990), no.~2, pp.~779--786.

\bibitem{dec-dioph}
rainmaker, \emph{Decidability of diophantine equation in a theory},
  MathOverflow, 2015, \url{http://mathoverflow.net/q/194491}.

\bibitem{sheph}
John~C. Shepherdson, \emph{A nonstandard model for a free variable fragment of
  number theory}, Bul\-le\-tin de l'Aca\-d{\'e}\-mie Po\-lo\-naise des
  \hbox{Sciences}, S{\'e}\-rie des \hbox{Sciences} Ma\-th{\'e}\-ma\-tiques,
  As\-tro\-no\-miques et Phy\-siques 12 (1964), no.~2, pp.~79--86.

\bibitem{wil:iop}
Alex~J. Wilkie, \emph{Some results and problems on weak systems of arithmetic},
  in: Logic {C}olloquium '77 (A.~Macintyre, ed.), North-Holland, 1978,
  pp.~285--296.

\end{thebibliography}
{\catcode`\/=13
  \gdef/{\string/\futurelet\nexttoken\finishslash}
  \gdef\finishslash{\ifx\nexttoken/\else\penalty\relpenalty\fi}
}
\providecommand\url{\begingroup\catcode`\~=12 \catcode`\/=13 \finishurl}
\def\finishurl#1{\texttt{#1}\endgroup}
\providecommand\bysame{\leavevmode\hbox to5em{\hrulefill}\thinspace}

\end{document}